\documentclass{article}
\usepackage[utf8]{inputenc}
\usepackage[english]{babel}
\usepackage{amsthm}
\usepackage{graphicx}
\usepackage{url}
\usepackage{hyperref}
\usepackage{amsmath}
\usepackage{xcolor}

\begin{document}

\title{Nested formulas for cosine and inverse cosine functions based on Vi\`{e}te's formula for $\pi$}

\author{Artur Kawalec}

\date{}
\maketitle

\begin{abstract}
In this article, we develop nested representations for cosine and inverse cosine functions, which is a generalization of Vi\`{e}te's formula for $\pi$. We explore a natural inverse relationship between these representations and develop numerical algorithms to compute them. Throughout this article, we perform numerical computation for various test cases, and demonstrate that these nested formulas are valid for complex arguments and a $k$th branch. We further extend the presented results to hyperbolic cosine and logarithm functions, and using additional trigonometric identities, we explore the sine and tangent functions and their inverses.
\end{abstract}

\section{Introduction}
In 1593s, Vi\`{e}te developed an infinite nested radical product formula for $\pi$ given as

\begin{equation}\label{eq:1}
\frac{2}{\pi}= \frac{\sqrt{2}}{2}\frac{\sqrt{2+\sqrt{2}}}{2}\frac{\sqrt{2+\sqrt{2+\sqrt{2}}}}{2}\cdots,
\end{equation}
which, historically, may have been the first infinite series for $\pi$. Vi\`{e}te derived this formula geometrically by means of approximating areas of an $n$-sided polygon inscribed in a circle, and obtaining a successive sequence of ratios of areas of $2n$-sided polygon to an $n$-sided polygon that is converging to an area of a circle, from which $2/\pi$ follows as a converging value as $n\to\infty$ [2]. Subsequently, Euler also developed a new proof of Vi\`{e}te's formula using trigonometric identities and infinite products, namely, Euler found a product identity

\begin{equation}\label{eq:1}
\frac{\sin(x)}{x}=\prod_{n=1}^{\infty}\cos\left(\frac{x}{2^n}\right)=\cos\left(\frac{x}{2}\right)\cos\left(\frac{x}{2^2}\right)\cos\left(\frac{x}{2^3}\right)\ldots,
\end{equation}
where by a repeated application of the half-angle identity for cosine

\begin{equation}\label{eq:1}
\cos\left(\frac{x}{2}\right)=\pm\sqrt{\frac{1-\cos(x)}{2}}
\end{equation}
the Vi\`{e}te's formula for $\pi$ is obtained by substituting $x=\pi/2$. Since then, there may have been perhaps dozens of variations of nested radical formulas published in the literature, and the main theme behind these formulas is a repeated application of the half or double angle identity for cosine function, which is a starting point of this article.

We will derive a nested formula for cosine function as
\begin{equation}\label{eq:1}
\cos(x) = \lim_{n\to\infty}\Bigg\{-1\ldots+2\ldots\Big(-1+2\Big(-1+2\Big(1-\frac{x^2}{2^{2n+1}}\Big)^2\Big)^2\Big)^{2\ldots}\Bigg\}
\end{equation}
which is an entire function, valid for complex $x$. In this notation, the dotted symbol $\ldots$ implies an $n$th nested recursion of the type $-1+2x^2$ until the inner term $1-x^2/2^{2n+1}$ is reached. Secondly, we will derive a nested formula for an inverse cosine function as

\begin{equation}\label{eq:1}
\cos^{-1}(y)=\lim_{n\to\infty} 2^{n}\sqrt{2\left(1- \sqrt{\frac{1}{2}\left(1+\sqrt{\frac{1}{2}\left(1+\sqrt{\frac{1}{2}\left(1+\sqrt{\frac{1}{2}\left(1+\ldots y\right)}\right)}\right)}\right)}\right)},
\end{equation}
that is also valid for a complex $y$. In this notation, the dotted symbol $\ldots$ implies an $n$th nested recursion of the type $\sqrt{\frac{1}{2}(1+y)}$ until the outer term $2^n\sqrt{2(1-y)}$ is reached. Hence, in this view, the Vi\`{e}te's formula is essentially a statement of $\cos^{-1}(0)=\pi/2$, but is originally presented as a reciprocal (1) and the square roots factored out. As it is seen, there is a natural inverse relationship between representations (4) and (5), as one can transform one form to the other by unraveling the outer terms into inner terms of its inverse. There is also closely related variation of (5) for an inverse cosine derived by Levin in [3] by means of a nested application of a certain hypergeometric function.

Throughout this article, we will also develop computer algorithms in Matlab software package to compute the presented formulas numerically. The simplicity of the code can easily be ported to any programming language. We will further extend the nested inverse cosine function (5) to be valid for a $k$th branch by means of a Gray code alternating sign function, which will be explored in more details. Although the main focus of this article is the cosine function and its inverse, we will further extend the main cosine formulas using basic trigonometric identities to other functions, such as the sine and the tangent, and their hyperbolic counterparts, as well as the nested logarithm function, and briefly mention the exponential function to complete the set of these nested formulas. Finally, we will summarize all the formulas in Appendix A.

\section{Nested formula for cosine function}
In this section, we will develop the nested formula for cosine (4). First, we re-write the half-angle identity for cosine function (3) as

\begin{equation}\label{eq:1}
\cos(x) = -1+2\cos^2\Big(\frac{x}{2}\Big),
\end{equation}
which can also be written as

\begin{equation}\label{eq:1}
\cos(x) = T_2\Big(\cos\Big(\frac{x}{2}\Big)\Big)
\end{equation}
for $n=1$ order, where the function $T_2(x)$ can be identified as a second order Chebyshev polynomial of the first kind as

\begin{equation}\label{eq:1}
T_2(x)= -1+2x^2,
\end{equation}
but we shall just refer to it as $T(x)$, since in this paper, we only work with the half-angle identity for cosine (6). More details about Chebyshev polynomials can be found in [1]. We then observe that if we recursively iterate $T(x)$ again, then we can recover the cosine function as

\begin{equation}\label{eq:1}
\cos(x) = T\Big(T\Big(\cos\Big(\frac{x}{2^2}\Big)\Big)\Big),
\end{equation}
which yields an expansion
\begin{equation}\label{eq:1}
\cos(x) = -1+2\Big(-1+2\cos^2\Big(\frac{x}{2^2}\Big)\Big)^2
\end{equation}
for $n=2$ order. If we repeat this process again, then we obtain

\begin{equation}\label{eq:1}
\cos(x) = T\Big(T\Big(T\Big(\cos\Big(\frac{x}{2^3}\Big)\Big)\Big)\Big),
\end{equation}
which yields an expansion
\begin{equation}\label{eq:1}
\cos(x) = -1+2\Big(-1+2\Big(-1+2\cos^2\Big(\frac{x}{2^3}\Big)\Big)^2\Big)^2
\end{equation}
for $n=3$ order. Hence, in general, we can write the cosine function as the $n$th order iterate of $T(x)$ as

\begin{equation}\label{eq:1}
\cos(x) = T^{\circ n}\Big(\cos\Big(\frac{x}{2^n}\Big)\Big),
\end{equation}
where the $n$th iteration on $T$ is defined as a recursive composition

\begin{equation}
T^{\circ k}(x) = T^{}(x)\circ T^{\circ k-1}(x),\quad k>0
\end{equation}
and
\begin{equation}
T^{\circ 0}(x) = f(x),\quad k=0,
\end{equation}
where $f(x)$ is the inner most function which is $\cos(\frac{x}{2^n})$, and it rapidly approaches $1$ as $n\to\infty$. We also know that the Maclaurin (Taylor) series for $\cos(x)$ about $x=0$ is
\begin{equation}\label{eq:1}
\cos(x) =1-\frac{x^2}{2!}+\frac{x^4}{4!}-\frac{x^6}{6!}+...,
\end{equation}
hence it suffices to use the first two terms as an approximation for the inner term as

\begin{equation}\label{eq:1}
\cos\Big(\frac{x}{2^n}\Big) \approx 1-\frac{x^2}{2^{2n+1}}.
\end{equation}
As a result, by substituting the inner term approximation (17) back to the infinite nested product (13) we obtain the main nested formula for cosine function

\begin{equation}\label{eq:1}
\cos(x) = \lim_{n\to \infty} T^{\circ n}\Big(1-\frac{x^2}{2^{2n+1}}\Big).
\end{equation}
We further observe that when the nested product formula (18) is expanded, it approaches the Maclaurin series for cosine (16). For example, for $n=0$ we start with inner approximation for cosine as

\begin{equation}\label{eq:1}
\cos(x) \approx 1-\frac{x^2}{2}.
\end{equation}
For $n=1$ we obtain an expansion

\begin{equation}\label{eq:18}
\begin{split}
\cos(x) & \approx -1+2\Big( 1-\frac{x^2}{2^3}\Big)^2 \\
        & \approx 1-\frac{1}{2}x^2+\frac{1}{32}x^4.
\end{split}
\end{equation}
For $n=2$ we have

\begin{equation}\label{eq:18}
\begin{split}
\cos(x) & \approx -1+2\Big(-1+2\Big( 1-\frac{x^2}{2^5}\Big)^2\Big)^2 \\
        & \approx 1-\frac{1}{2}x^2+\frac{5}{128}x^4-\frac{1}{1024}x^6+\frac{1}{131072}x^8.
\end{split}
\end{equation}
And, for $n=3$ we have
\begin{equation}\label{eq:18}
\begin{split}
\cos(x) & \approx -1+2\Big(-1+2\Big(-1+2\Big(1-\frac{x^2}{2^7}\Big)^2\Big)^2\Big)^2 \\
        & \approx 1-\frac{1}{2}x^2+\frac{21}{512}x^4-\frac{21}{16384}x^6+\frac{165}{8388608}x^8-...,
\end{split}
\end{equation}
where higher order terms from $10$ to $16$ have been omitted. In fact, $2^n+1$ terms of approximated Maclaurin series are generated, which approach the Maclaurin series as $n\to\infty$. One major advantage of this iterative process is that the factorial for large $n$ doesn't have to be computed if one were to use the Maclaurin series to compute the cosine function.

In the next example, we numerically compute $\cos(\pi/3)$ by application of equation (18) for $n=4$. First we begin to compute the inner term $y_0$, and then recursively apply the Chebyshev polynomial $T(x)$ until we reach the outer term as

\begin{equation}\label{eq:18}
\begin{split}
y_0     & = 1-(\pi/3)^2/2^9 =   0.997858158767125, \\
y_1        & = -1+2y_0^2 =     0.991441810036233, \\
y_2        & = -1+2y_1^2 =   0.965913725375842, \\
y_3        & = -1+2y_2^2 =    0.865978649738875, \\
y_4        & = -1+2y_3^2 =    0.499838043607131, \\
\end{split}
\end{equation}
hence we obtain the final result as

\begin{equation}\label{eq:4}
\cos(\pi/3) \approx  0.499838043607131,
\end{equation}
where we obtain accuracy to within three decimal places. Such accuracy could be reached with $4$ terms of Maclaurin series. In next example, we repeat the computation for $n=10$ case as

\begin{equation}\label{eq:18}
\begin{split}
y_0     & = 1-(\pi/3)^2/2^{21} =   0.999999477089543,\\
y_1      &  = -1+2y_0^2  =      0.999997908358718,\\
y_2        & = -1+2y_1^2 =      0.999991633443621,\\
y_3        & = -1+2y_2^2 =      0.999966533914483,\\
y_4        & = -1+2y_3^2 =      0.999866137897891,\\
y_5        & = -1+2y_4^2 =      0.999464587429687,\\
y_6        & = -1+2y_5^2 =      0.997858923051989,\\
y_7        & = -1+2y_6^2 =      0.991444860628951,\\
y_8        & = -1+2y_7^2 =      0.965925823335118,\\
y_9        & = -1+2y_8^2 =      0.866025392371252,\\
y_{10}        & = -1+2y_9^2 =   0.499999960463562,\\
\end{split}
\end{equation}
where we obtain accuracy to within $7$ decimal places. The number of term of Maclaurin series to reach such accuracy is greater than $7$, which requires computing factorial of $14$ which is on the order of $10^{11}$, while $2^{21}$ is on the order of $10^6$. 

We can further improve convergence by approximating the inner cosine terms with more terms of Maclaurin series, for example, with $4$ terms we have an approximation to cosine as

\begin{equation}\label{eq:1}
\cos_4\Big(\frac{x}{2^n}\Big) \approx 1-\frac{x^2}{2^{2n}2!}+\frac{x^4}{2^{4n}4!}-\frac{x^6}{2^{6n}6!}.
\end{equation}
Hence, we again recompute the cosine as

\begin{equation}\label{eq:18}
\begin{split}
y_0     & = \cos_4((\pi/3)/2^4) =   0.997858923238595, \\
y_1        & = -1+2y_0^2 =      0.991444861373777, \\
y_2        & = -1+2y_1^2 =      0.965925826288936,\\
y_3        & = -1+2y_2^2 =      0.866025403783926,\\
y_4        & = -1+2y_3^2 =      0.499999999998225, \\
\end{split}
\end{equation}
where now we reached accuracy to within $12$ decimal places. In essence, this method is similar to a convergence acceleration algorithm.  Also, since it converges to Maclaurin series for the cosine function, it is therefore an entire function valid in all complex plane.

Next, we illustrate a practical algorithm developed in the Matlab software package, and the code can be easily implemented in any programming language. We define a function $\text{cos\_fx(x,k)}$:
\begin{verbatim}
1:   % This function evaluates cos(x)
2:   % x is the input argument (real or complex)
3:   % k=2^n (must be power of 2) where n is the order of nested product
4:   %
5:   function [y] = cos_fx(x,k)
6:
7:       if k==1
8:           y = 1-x^2/2;                % Inner term approximation to cosine
9:       else
10:          y = -1+2*cos_fx(x/2,k/2)^2; % Recursively call Chebyshev polynomial of cos_fx
11:      end
12:
13:  end
\end{verbatim}
The code is similarly implemented as a factorial function recursively calling itself. We call this function for $n=4$ as
\begin{verbatim}
cos_fx(pi/3,2^4) = 0.499838043607131
\end{verbatim}
where $k=2^n$ and must be a power of $2$, and we thus repeat the result of (23). And similarly, we can improve the convergence by adding more terms of the initial approximation to cosine as

\begin{verbatim}
1:  % This function evaluates cos(x)
2:  % x is the input argument (real or complex)
3:  % k=2^n (must be power of 2) where n is the order of nested product
4:  %
5:  function [y] = cos_fx4(x,n)
6:
7:      if n==1
8:          y = 1-x^2/2+x^4/24-x^6/720; % Inner term approximation to cosine
9:      else
10:          y = -1+2*cos_fx4(x/2,n/2)^2; % Chebyshev polynomial of cos_fx4
11:     end
12:
13: end
\end{verbatim}
which for $n=4$ results in
\begin{verbatim}
cos_fx4(pi/3,2^4)=  0.499999999998225
\end{verbatim}
thus repeating the result of (27).

\section{The nested inverse cosine}

The nested product structure for the cosine function permits reversing the recursive process, resulting in a formula for an inverse cosine function. To illustrate this, we consider the inner cosine approximation (16) as

\begin{equation}\label{eq:1}
y = \cos(x) \approx 1-\frac{x^2}{2},
\end{equation}
then by solving for $x$, we obtain an initial approximation to inverse cosine function $x =\cos^{-1}(y)$ as

\begin{equation}\label{eq:1}
g(y)= \sqrt{2(1-y)},
\end{equation}
which we consider for $n=0$ order approximation. We also find an inverse of the Chebyshev polynomial $T(x)$ as

\begin{equation}\label{eq:1}
T^{-1}(y)=\sqrt{\frac{1+y}{2}}.
\end{equation}
If we consider the next higher order approximation (20) for $n=1$, then we have

\begin{equation}\label{eq:1}
y = \cos(x)\approx -1+2\Big( 1-\frac{x^2}{2^3}\Big)^2,
\end{equation}
and by solving for $x$ results in
\begin{equation}\label{eq:1}
x = \cos^{-1}(y)\approx 2\sqrt{2\left(1-\sqrt{\frac{1}{2}(1+y)}\right)}.
\end{equation}
If we consider the next higher order approximation (21) for $n=2$, then we have
\begin{equation}\label{eq:1}
y = \cos(x) \approx -1+2\Big(-1+2\Big( 1-\frac{x^2}{2^5}\Big)^2\Big)^2
\end{equation}
and by solving for $x$ results in
\begin{equation}\label{eq:1}
x = \cos^{-1}(y)\approx 2^2\sqrt{2\left(1-\sqrt{\frac{1}{2}\left(1+\sqrt{\frac{1}{2}(1+y)}\right)}\right)}.
\end{equation}
If we consider the next higher order approximation (22) for $n=3$, then we have
\begin{equation}\label{eq:1}
y = \cos(x) \approx -1+2\Big(-1+2\Big(-1+2\Big(1-\frac{x^2}{2^7}\Big)^2\Big)^2\Big)^2
\end{equation}
and by solving for $x$ results in
\begin{equation}\label{eq:1}
x = \cos^{-1}(y)\approx 2^3\sqrt{2\left(1-\sqrt{\frac{1}{2}\left(1+\sqrt{\frac{1}{2}\left(1+\sqrt{\frac{1}{2}(1+y)}\right)}\right)}\right)}.
\end{equation}
Hence, continuing on for an $n$th iteration, it is seen that by solving for $x$ unravels the nested formula for cosine which directly results in a nested formula for inverse cosine function as
\begin{equation}\label{eq:1}
\cos^{-1}(y)=\lim_{n\to\infty} 2^{n}\sqrt{2\left(1-\sqrt{\frac{1}{2}\left(1+\sqrt{\frac{1}{2}\left(1+\sqrt{\frac{1}{2}\left(1+\sqrt{\frac{1}{2}(1+\ldots y)}\right)}\right)}\right)}\right)},
\end{equation}
which also can be written as
\begin{equation}\label{eq:1}
\cos^{-1}(y) = \lim_{n\to \infty} 2^n g(T^{-1\circ n}(y)),
\end{equation}
where $g(y)=\sqrt{2(1-y)}$ is the outer function evaluated only once at the end, and the inverse function $T(y)^{-1}$ is recursively iterated an $n$ number of times. Now it is clear that Vi\`{e}te's formula is essentially $\cos^{-1}(0)=\pi/2$, but was originally presented as a reciprocal (1). Next, we numerically compute it for $n=4$ as

\begin{equation}\label{eq:18}
\begin{split}
x_0     & = \sqrt{\frac{1}{2}} =  0.707106781186548, \\
x_1        & = \sqrt{\frac{1+x_0}{2}} =      0.923879532511287, \\
x_2        & = \sqrt{\frac{1+x_1}{2}} =      0.980785280403230,\\
x_3        & = \sqrt{\frac{1+x_2}{2}} =       0.995184726672197,\\
x_4        & = 2^4\sqrt{2(1-x_3)} =      1.570165578477370, \\
\end{split}
\end{equation}
hence, the result is converging to $\pi/2$. As another example, we compute $\cos^{-1}(0.5)$ as

\begin{equation}\label{eq:18}
\begin{split}
x_0     & = \sqrt{\frac{1+0.5}{2}} =  0.866025403784439, \\
x_1        & = \sqrt{\frac{1+x_0}{2}} =      0.965925826289068, \\
x_2        & = \sqrt{\frac{1+x_1}{2}} =      0.991444861373810,\\
x_3        & = \sqrt{\frac{1+x_2}{2}} =      0.997858923238603,\\
x_4        & = 2^4\sqrt{2(1-x_3)} =      1.047010650296843, \\
\end{split}
\end{equation}
hence, the result is converging to $\pi/3$. As a result, we implement a simple function in Matlab $\text{acos\_fx(y,n)}$ as

\begin{verbatim}
1:   % This function evaluates acos(y)
2:   % y is the input argument (real or complex)
3:   % n is the order of nested product
4:   %
5:   function [x] = acos_fx(y,n)
6:
7:      for i = 1:n
8:
9:          y = sqrt((y+1)/2)  % Iterate inverse f(y) n number of times
10:
11:     end
12:
13:     x = 2^n*sqrt(2*(1-y))  % Outer term inverse cosine approximation g(y)
14:
15:  end
\end{verbatim}
We call this function for above values for $n=10$ to obtain

\begin{verbatim}
acos_fx(0,10) =  1.570796172805538
\end{verbatim}
and
\begin{verbatim}
acos_fx(0.5,10) = 1.047197505529385
\end{verbatim}
where we obtain accuracy to $7$ decimal places.

Generally, the cosine function maps a real number interval $(-\infty,+\infty)$ onto a unit interval $[-1,1]$, while the inverse cosine function is a multi-valued function that maps the unit interval $[-1,1]$ onto a real interval $(-\infty,+\infty)$ in small sections, or branches, which are enumerated by an integer $k$. Hence, if $y$ is on a unit interval $[-1,1]$, then for a $k$th branch, the mapping is

\begin{equation}\label{eq:1}
\cos^{-1}(y)= x + 2\pi k.
\end{equation}
For a principal branch $k=0$, the mapping is $[-1,1]\to [\pi,0]$. For $k=1$, the mapping is $[-1,1]\to [3\pi,2\pi]$ and for $k=2$ the mapping is $[-1,1]\to [5\pi,4\pi]$, and so on.
Also, the inverse cosine function is well-defined for a complex argument $z$ as

\begin{equation}\label{eq:1}
\cos^{-1}(z)=\frac{\pi}{2}+i\log{(iz+\sqrt{1-z^2})}
\end{equation}
having a branch cut on real axis in an interval $(-\infty,-1)$ and $(1,\infty)$. We find that the formula for inverse cosine (37) converges to (42). For example, we take

\begin{equation}\label{eq:1}
\cos^{-1}(2)= 1.316957896924817i
\end{equation}
while the formula (37) for $n=10$ gives
\begin{verbatim}
acos_fx(2,10) =  1.316956719106592i
\end{verbatim}
As another example, we take

\begin{equation}\label{eq:1}
\cos^{-1}(2+3i)=  1.000143542473797 - 1.983387029916535i
\end{equation}
while the formula (37) for $n=10$ gives
\begin{verbatim}
acos_fx(2,10) =   1.000533856110922 - 1.982613299971578i
\end{verbatim}
We performed higher precision computation and it clearly converged to (37) for any value we tried.

We further observe that since the square root function is multi-valued, generating a $\pm $ sign, hence the infinite nested radical product actually generates an infinite number of formulas with varying signs before the square root. We begin by changing the inner most sign as indicated in red color, and observe the convergence of

\begin{equation}\label{eq:1}
\cos^{-1}(y)=\lim_{n\to \infty} 2^{n}\sqrt{2\left(1-\sqrt{\frac{1}{2}\left(1+\sqrt{\frac{1}{2}\left(1+\sqrt{\frac{1}{2}\left(1 \mathbin{\color{red}{-}} \sqrt{\frac{1}{2}(1+\ldots y)}\right)}\right)}\right)}\right)},
\end{equation}
which results in for $y=0$
\begin{equation}\label{eq:1}
x = \cos^{-1}(0)\approx  4.712384822113464
\end{equation}
for $n=10$, where it is seen converging to $3\pi/2$. And similarly if we change the second inner most sign as

\begin{equation}\label{eq:1}
\cos^{-1}(y)= \lim_{n\to \infty} 2^{n}\sqrt{2\left(1-\sqrt{\frac{1}{2}\left(1+\sqrt{\frac{1}{2}\left(1\mathbin{\color{red}{-}}\sqrt{\frac{1}{2}\left(1+\sqrt{\frac{1}{2}(1+\ldots y)}\right)}\right)}\right)}\right)}
\end{equation}
for $n=10$ result in

\begin{equation}\label{eq:1}
x = \cos^{-1}(0)\approx 10.995521462263701,
\end{equation}
which is seen converging to $7\pi/2$. And similarly if we change two inner most signs as

\begin{equation}\label{eq:1}
\cos^{-1}(y)= \lim_{n\to \infty} 2^{n}\sqrt{2\left(1-\sqrt{\frac{1}{2}\left(1+\sqrt{\frac{1}{2}\left(1\mathbin{\color{red}{-}}\sqrt{\frac{1}{2}\left(1\mathbin{\color{red}{-}}\sqrt{\frac{1}{2}(1+\ldots y)}\right)}\right)}\right)}\right)}
\end{equation}
which for $n=10$ results in

\begin{equation}\label{eq:1}
x = \cos^{-1}(0)\approx 7.853962382754363,
\end{equation}
which is seen approaching $5\pi/2$, thus the signs of the square roots produces a sort of binary weighted code. As a result, we map this combination in Table 1 for a few more values by binary weighting the sign of the square root such that the LSB (Least Significant Bit) is a sign of the inner most square root and working outwards to the next sign. For the first four inner most square roots of (37) would be represented as $++++$, for equation (45) it would be $+++-$, for equation (47) it would be $++-+$, and for equation (49) it would be $++--$. From Table 1, it is seen that the pattern of signs resembles that of a Gray code if one takes $+$ to be $0$ and $-$ to be $1$. The Gray code has many applications in digital communication and computer science, such as error correction and data transmission. The $4$-bit Gray code for the previous sequence would be $0000$, $0001$, $0011$ and $0010$ for equation (37,45,49,47) respectively. As a result, we define a Gray code alternating sign function $G_k$ as
\begin{table}[ht]
\caption{Evaluation of Equation (37) for different square root sign} 
\centering 
\begin{tabular}{c c c} 
\hline\hline 
Square root sign & $\cos^{-1}(0)$ Eq(47) for $n=10$ & Converging value\\ [0.5ex] 
\hline 
$++++$  & 1.570796172805538  & $\pi/2$   \\
$+++-$  & 4.712384822113464  & $3\pi/2$  \\
$++--$  & 7.853962382754363  & $5\pi/2$  \\
$++-+$  & 10.995521462263701 & $7\pi/2$ \\
$+--+$  & 14.137054668250807 & $9\pi/2$ \\
$+---$  & 17.278554608371589 & $11\pi/2$ \\
$+-+-$  & 20.420013890390901 & $13\pi/2$\\
$+-++$  & 23.561425122147117 & $15\pi/2$
\\ [1ex] 
\hline 
\end{tabular}
\label{table:nonlin} 
\end{table}

\begin{equation}
G_k(m)= \left \{
\begin{aligned}
&1, &&\text{if}\ \text{$m$th binary digit of Gray code of $k$ is zero}\\
-&1,&& \text{if}\ \text{$m$th binary digit of Gray code of $k$ is one}.
\end{aligned} \right.
\end{equation}
to generate the $\pm 1$ in the $m$th nested radical. The main formula for a $k$th branch becomes now

\begin{equation}\label{eq:1}
\begin{aligned}
&\cos^{-1}(y)= \\
& \pm 2^{n}\sqrt{2\left(1-G_k(n)\sqrt{\frac{1}{2}\left(1+ G_k(3)\sqrt{\frac{1}{2}\left(1+G_k(2)\sqrt{\frac{1}{2}\left(1+G_k(1)\sqrt{\frac{1}{2}(1+\ldots y)}\right)}\right)}\right)}\right)}
\end{aligned}
\end{equation}
and in terms of an $n$th recursive composition can be written
\begin{equation}\label{eq:1}
T_{G_k}^{-1\circ (n)}(y) = G_k(n-1)T^{-1\circ (n-1)} (y)
\end{equation}
so that
\begin{equation}\label{eq:1}
\cos^{-1}(y) = \lim_{n\to \infty} \pm 2^n g(T_{G_k}^{-1\circ n}(y))
\end{equation}
is a full representation of an inverse cosine function for a $k$th branch. The principal branch $k=0$ has all square roots with positive sign except the outer function g(y). The plus and minus sign is due to the outer square root, where the positive sign is for branches for positive values of $k$, and the negative sign is for branches for negative values of $k$. We next will modify the algorithm for the inverse cosine valid for a $k$ branch as

\begin{verbatim}
1:   % This function evaluates acos(y)
2:   % y is the input argument (real or complex)
3:   % k is the kth branch k>=0
4:   % n is the order of nested product
5:   %
6:   function [x] = acos_fx_k(y,k,n)
7:
8:      Gk = Gk_vec(k,n);  % Load Gray code sequence
9:
10:     for i = 1:n
11:
12:         y = Gk(i)*sqrt((y+1)/2);   % Iterate Gk * inverse T(y)
13:
14:     end
15:
16:     x = 2^n*sqrt(2*(1-y));   % Apply outer inverse cosine approximation
17:
18:  end
\end{verbatim}
And the Gray code alternating sign function is $G_{k}(m)$ is given as:

\begin{verbatim}
1:   % This function generates Gray code alternating sign vector
2:   % m is the index of Gray code sequence and must be positive integer
3:   % n is the order of nested product and length of Gray code vector
4:   %
5:  function [Gk] = Gk_vec(m,n)
6:
7:       bin = dec2bin(m,n);     % Convert integer value m to binary string of length n
8:       gcode = ones(1,n);      % Initialize vector with '1' to store Gk values
9:       gcode_refl = ones(1,n); % Initialize vector with '1' to store reflected Gk values
10:
11:      for i = 2:n
12:
13:        bit1 = str2double(bin(i-1)); % i-1 binary bit1 string converted to double
14:        bit2 = str2double(bin(i));   % i binary bit2 string converted to double
15:
16:        % Apply xor operation on bit1 and bit2
17:        c = xor(bit1,bit2);
18:
19:        % Alternating sign value assignment
20:        if c == 0
21:            gcode(i) = 1;
22:        elseif c ==1
23:            gcode(i) = -1;
24:        end
25:
26:     end
27:
28:     % Reflect alternating sign values so that MSB corresponds to LSB
29:     for i = 1:n
30:         gcode_refl(n-i+1) = gcode(i);
31:     end
32:
33:     % Return Gk
33:     Gk = gcode_refl;
34:
35: end
\end{verbatim}
where we utilized some unique functions in Matlab, such as `dec2bin' function to convert a decimal value into a vector filled with corresponding binary elements, and the `ones' function allocates a vector filled with values of $1$. The Gray code sequence is generated by a XOR operation on two successive bits.

When we test this function, we quickly reproduce all the values in Table 1, and next, we try computing the inverse cosine for the $100$th branch for $n=25$. The Gray code alternating sign function results in a sequence

\begin{equation}\label{eq:1}
G_{100} = \{1,1,1,1,1,1,1,1,1,1,1,1,1,1,1,1,1,1,-1, 1,-1,1,-1,-1, 1\}
\end{equation}
where we use this relation to extract the $k$th branch
\begin{equation}\label{eq:1}
 k = \frac{\cos^{-1}(0)}{\pi}-\frac{1}{2}
\end{equation}
by dividing by $\pi$ and subtracting $\frac{1}{2}$ to obtain
\begin{verbatim}
acos_fx_k(0,100,25)/pi-0.5 =  100.000014188946
\end{verbatim}
accurate to four decimal places. And similarly, for the millionth branch $k=10^6$ for $n=25$ we have

\begin{equation}\label{eq:1}
G_{10^6} = \{1,1,1,1,1,-1,1,1,1,-1,-1,-1,1,1,1,-1,-1,1,-1,-1,1,1,1,1,1\}
\end{equation}
and we obtain
\begin{verbatim}
acos_fx_k(0,10^6,25)/pi-0.5 =  999634.790737135
\end{verbatim}
while for $n=30$ we obtain

\begin{verbatim}
acos_fx_k(0,10^6,30)/pi-0.5 =  999999.643311845,
\end{verbatim}
an accuracy to almost $1$ decimal place with just $30$ nested radicals for a millionth branch. We also performed a sweep of the branch $k$ from $0$ to $10^6$ in $1$ step and plotted the extracted branches  in Fig. 1 to verify the validity of the Gray code sequence $G_k$. Indeed, we see a perfect reproduction up to a millionth branch, only with $n=25$ nested radicals. Also, the $kth$ branch  works well for a complex argument, we compute the millionth branch again for $\cos^{-1}(2+3i)$ with $n=30$ and obtain
\begin{verbatim}
real(acos_fx_k(2+3*i,10^6,30)/pi) =  999999.961666679
\end{verbatim}
where we roughly obtain the millionth extracted branch, since we don't know the exact multiple of $\pi$ to subtract as in (55) case.

\begin{figure}[hbt!]
  \centering
  \includegraphics[width=100mm]{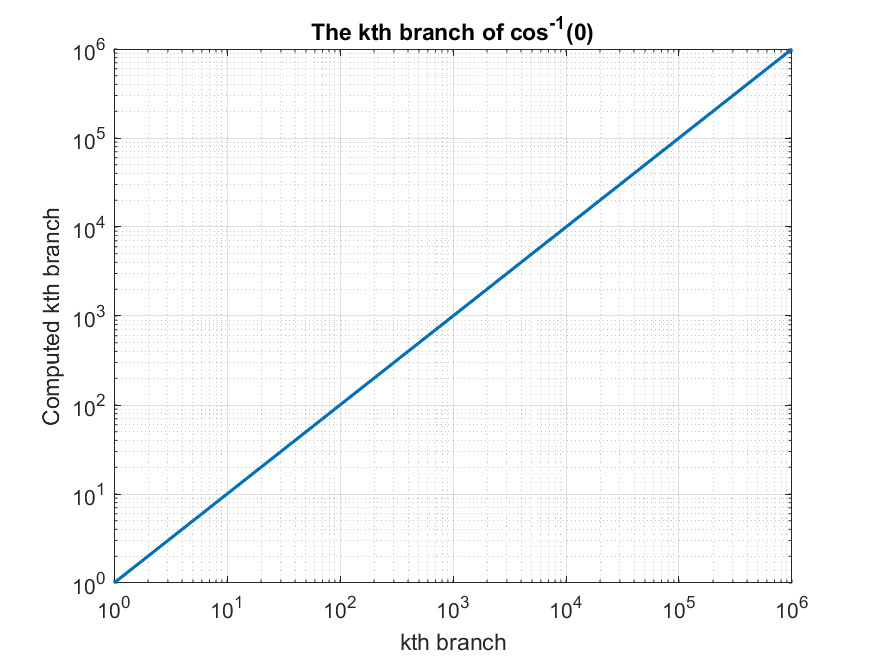}\\
  \caption{Extracted $k$th branch of $\cos^{-1}(0)$ by equation (52) from $k=0$ to $k=10^6$ in unit steps}\label{1}
\end{figure}

There is also one special case for $\pm 1$ arguments. Namely, for $\cos^{-1}(1)$ and an even branch $k$, there is also the same $k-1$ branch. And similarly, for $\cos^{-1}(-1)$ and an odd branch, there is also same $k-1$ branch. We summarize this in Table 2 by computing $k=\cos^{-1}(\pm 1)/\pi$ and note how the double branch arise for even and odd branches.  From equation (37) it cancels or maximizes the inner most radical term and the effect propagates up the radical chain, since the LSB bit of the Gray code doesn't change for one increment of $k$. 

\begin{table}[hbt!]
\caption{The $k$th branch of $\cos^{-1}(-1)$ and $\cos^{-1}(1)$ by equation (52)} 
\centering 
\begin{tabular}{c c c} 
\hline\hline 
kth branch & $\cos^{-1}(1)/\pi$  & $\cos^{-1}(-1)/\pi$\\ [0.5ex] 
\hline 
$0$  & 0 & 0.999999607815114  \\
$1$  & 1.99999686254052 & 0.999999607815114 \\
$2$  & 1.99999686254052 & 2.99998941107727 \\
$3$  & 3.99997490034631 & 2.99998941107727 \\
$4$  & 3.99997490034631 & 4.99995097728807 \\
$5$  & 5.99991528886506 & 4.99995097728807 \\
$6$  & 5.99991528886506 & 6.99986548205935 \\
$7$  & 7.99979920389600 & 6.99986548205935 \\
$8$  & 7.99979920389600 & 8.99971410143166 \\
$9$  & 9.99960782177148 & 8.99971410143166 \\
$10$ & 9.99960782177148 & 10.9994780120676
\\ [1ex] 
\hline 
\end{tabular}
\label{table:nonlin} 
\end{table}

\section{Hyperbolic Cosine}
The hyperbolic cosine function is defined as
\begin{equation}\label{eq:1}
\cosh(x) = \cos(ix)
\end{equation}
and since the nested product formula (4) approaches the Maclaurin series for cosine function, it is therefore an entire function valid in all complex plane, as a result, we then have
\begin{equation}\label{eq:1}
\cosh(x) = \lim_{n\to \infty} T^{\circ n}\Big(1+\frac{x^2}{2^{2n+1}}\Big),
\end{equation}
which yields an expansion of

\begin{equation}\label{eq:1}
\cosh(x) = \lim_{n\to\infty}\Bigg\{-1\ldots+2\ldots\Big(-1+2\Big(-1+2\Big(1+\frac{x^2}{2^{2n+1}}\Big)^2\Big)^2\Big)^{2\ldots}\Bigg\}.
\end{equation}
And similarly as before, by solving for $x$ we obtain an inverse function as

\begin{equation}\label{eq:1}
\cosh^{-1}(y)=\lim_{n\to\infty} 2^{n}\sqrt{2\left(-1+\sqrt{\frac{1}{2}\left(1+\sqrt{\frac{1}{2}\left(1+\sqrt{\frac{1}{2}\left(1+\sqrt{\frac{1}{2}(1+\ldots y}\right)}\right)}\right)}\right)},
\end{equation}
where only the outer function approximation to inverse cosine $g(y)$ changed to an approximation to an inverse hyperbolic cosine as $h(y)=\sqrt{2(-1+y)}$. Also, in terms of the branches, we have

\begin{equation}\label{eq:1}
\cosh^{-1}(y)=\lim_{n\to\infty} \pm 2^n h(T_{G_k}^{-1\circ n}(y))
\end{equation}
where $G_k$ is the same Gray code ordering of the square root signs establishes the branches, such that in $G_0$ all square root signs are positive, in $G_1$ only the first inner most sign is negative, in $G_2$ only the first two inner most signs are negative and in $G_3$ only the second inner most sign is negative and so on. Also, the inverse hyperbolic cosine function is defined for complex argument $z$ as

\begin{equation}\label{eq:1}
\cosh^{-1}(z) = \log(z+\sqrt{z^2-1})
\end{equation}
having a branch cut on real axis in an interval $(-\infty,1)$. To test this function we compute

\begin{equation}\label{eq:1}
\cosh^{-1}(2) =  1.31695789692482
\end{equation}
while the computer function for $n=10$ gives
\begin{verbatim}
acosh_fx(2,10) =   1.31695798760619
\end{verbatim}
As another example, if we take

\begin{equation}\label{eq:1}
\cosh^{-1}(-2+3i) =  1.98338702991654 - 2.141449111116i
\end{equation}
while the computer function for $n=10$ gives
\begin{verbatim}
acosh_fx(-2+3i,10) = 1.98338625571006 - 2.14144972510396i
\end{verbatim}

Finally, we give a slightly modified version of the code to compute the hyperbolic cosine function, where we use $4$ terms to approximate the inner $\cosh(x)$ approximation
\begin{verbatim}
1:   % This function evaluates cosh(x)
2:   % x is the input argument (real or complex)
3:   % k=2^n (must be power of 2) where n is the order of nested product
4:   %
5:   function [y] = cosh_fx4(x,n)
6:
7:      if n==1
8:          y = 1+x^2/2+x^4/24+x^6/720;  % Inner term approximation to hyperbolic cosine
9:      else
10:         y = -1+2*cosh_fx(x/2,n/2)^2; % Chebyshev polynomial of cosh_fx
11:     end
12:
13:  end
\end{verbatim}
and the inverse function as

\begin{verbatim}
1:   % This function evaluates acosh(y)
2:   % y is the input argument (real or complex)
3:   % n is the order of nested product
4:   %
5:   function [x] = acosh_fx(y,n)
6:
7:      for i = 1:n
8:
9:          y = sqrt((y+1)/2)  % Iterate inverse f(y) n number of times
10:
11:     end
12:
13:     x = 2^n*sqrt(2*(-1+y))  % Outer term inverse acosh approximation h(y)
14:
15:  end.
\end{verbatim}
But, if one wishes to add the branches, then the $G_k$ must be included simiarly as for the inverse cosine function in the previous section.

\section{The logarithm function}
By observing the definition of inverse hyperbolic cosine (63) as

\begin{equation}\label{eq:1}
\cosh^{-1}(x) = \log(x+\sqrt{x^2-1})
\end{equation}
we write the inner log term as
\begin{equation}\label{eq:1}
y = x+\sqrt{x^2-1}
\end{equation}
from which solving for $x$ yields
\begin{equation}\label{eq:1}
x = \frac{y^2+1}{2y}.
\end{equation}
This now gives a formula for log in terms of inverse hyperbolic cosine as
\begin{equation}\label{eq:1}
\log(y) = \cosh^{-1}\left(\frac{y+y^{-1}}{2}\right),
\end{equation}
and by using the nested radical representation for inverse hyperbolic cosine (61), the logarithm formula becomes

\begin{equation}\label{eq:1}
\log(y) =\lim_{n\to\infty} 2^{n}\sqrt{2\left(-1+\sqrt{\frac{1}{2}\left(1+ \sqrt{\frac{1}{2}\left(1+\sqrt{\frac{1}{2}\left(1+\sqrt{\frac{1}{2}\left(1+\ldots \frac{y+y^{-1}}{2}\right)}\right)}\right)}\right)}\right)},
\end{equation}
where, again as before, the inner inverse Chebyshev function $T^{-1}(x)$ is iterated $n$ number of times as

\begin{equation}\label{eq:1}
\begin{aligned}
\log(y) & = \lim_{n\to \infty} 2^n h\Big(T^{-1}\enskip \Big(T^{-1}\Big(T^{-1}\Big(T^{-1}\Big(\ldots\frac{y+y^{-1}}{2}\Big)\Big)\Big)\Big)\Big) \\
        & = \lim_{n\to \infty} 2^n h\Big(T^{-1\circ n}\Big(\frac{y+y^{-1}}{2}\Big)\Big),
\end{aligned}
\end{equation}
and $h(y)$ is the outer function evaluated only once in the end. If one employs the Gray code alternating sign sequence, then one could recover all the branches of the logarithm. Next, we investigate this formula for $\log(2)$ and find that inner term becomes

\begin{equation}\label{eq:1}
\log(2) =\lim_{n\to\infty} 2^{n}\sqrt{2\left(-1+\sqrt{\frac{1}{2}\left(1+\sqrt{\frac{1}{2}\left(1+\sqrt{\frac{1}{2}\left(1+\sqrt{\frac{1}{2}\left(1+\ldots \frac{5}{4}\right)}\right)}\right)}\right)}\right)}.
\end{equation}
A computation for $n=10$ nested radicals yields
\begin{verbatim}
log_fx(2,10) =   0.693147193652913
\end{verbatim}
with an accuracy to $6$ decimal places. We coded the function $\text{log\_fx}$ which is just a simple modification  of $\text{acosh\_fx}$ shown in the previous section. For $\log(-1)$, interestingly, the inner radical term becomes zero, and we are left with a nested radical formula for $\cosh^{-1}(-1)=\pi i$, and when we reach the outer radical term in which there is a square root of $-1$ to give the imaginary part. To verify this, we compute
\begin{verbatim}
log_fx(-1,10) =   3.14159142150464i
\end{verbatim}
for $n=10$ with an accuracy to $5$ decimal places. Next, we demonstrate how $\log(1)=0$, where it is seen the inner most term becomes
\begin{equation}\label{eq:1}
\sqrt{\frac{1}{2}\left(1+\frac{1+1}{2}\right)}=1
\end{equation}
the next outer radical term is
\begin{equation}\label{eq:1}
\sqrt{\frac{1}{2}\left(1+\frac{1}{2}2\right)}=1
\end{equation}
and this process is propagated up the radical chain like a domino effect, until we reach the outer most radical term where it is terminated with $-1+1=0$.

Finally, we note that this formula can handle very large input argument. For example, we compute $\log(10^{100})= 230.258509299405$ also for $n=10$, and obtain

\begin{verbatim}
log_fx(1e100,10) =  230.743921174299.
\end{verbatim}
Although, in practice, one would use the rule $\log(a^x)=x\log(a)$ to simplify computation of logarithms for large argument, however it is seen that this nested radical formula can cut through large input argument very effectively.

\section{The sine, tangent and exponential functions}
The the half-angle identify for sine
\begin{equation}\label{eq:1}
\sin(x) = 2\sin\left(\frac{x}{2}\right)\cos\left(\frac{x}{2}\right)
\end{equation}
relates the sine function to the product of sine with cosine functions, thus preventing to recursively iterate the sine function similarly as in (4). However, what one might could do is to write 

\begin{equation}\label{eq:1}
\sin(x) = \cos\left(x-\frac{\pi}{2}\right)
\end{equation}
and use the nested cosine formula. Another possibility is apply the identity $\cos^2(x)+\sin^2(x)=1$ so that

\begin{equation}\label{eq:1}
\sin(x) = \sqrt{1-\cos^2(x)},
\end{equation}
which results in a slight modification

\begin{equation}\label{eq:1}
\sin(x) = \sqrt{1-\lim_{n\to\infty}\Bigg\{-1\ldots+2\ldots\Big(-1+2\Big(-1+2\Big(1-\frac{x^2}{2^{2n+1}}\Big)^2\Big)^2\Big)^{2\ldots}\Bigg\}^2},
\end{equation}
where the square root term is the outer most term. And for the inverse sine, we use the following identity

\begin{equation}\label{eq:1}
\sin^{-1}(y) = \cos^{-1}(\sqrt{1-y^2}),
\end{equation}
and by substituting to the nested inverse cosine function we obtain

\begin{equation}\label{eq:1}
\sin^{-1}(y)= \lim_{n\to\infty} 2^{n}\sqrt{2\left(1-\sqrt{\frac{1}{2}\left(1+ \sqrt{\frac{1}{2}\left(1+\sqrt{\frac{1}{2}\left(1+\sqrt{\frac{1}{2}\left(1+\ldots \sqrt{1-y^2}\right)}\right)}\right)}\right)}\right)},
\end{equation}
where the inner most term has the square root. In terms of the recursive composition is represented as

\begin{equation}\label{eq:1}
\sin^{-1}(y) = \lim_{n\to \infty}  2^n g(T^{-1\circ n}(\sqrt{1-y^2})).
\end{equation}
The sine functions are basically coupled to the cosine function by equation (77). A similar situation occurs for the tangent function

\begin{equation}\label{eq:1}
\tan(x) = \frac{2\tan(\frac{1}{2}x)}{1-\tan^2{(\frac{1}{2}x})}= \frac{2}{\tan(\frac{1}{2}x)-\frac{1}{\tan(\frac{1}{2}x)}}
\end{equation}
where the iterated half-angle identity generates a nested tree that would grow similarly as a nested continued fraction. So at this point we just refer back to the tangent function as a combination of sine and cosine as

\begin{equation}\label{eq:1}
\tan(x) = \frac{\sin(x)}{\cos(x)}= \frac{\sqrt{1-\cos^2(x)}}{\cos(x)}= \sqrt{\frac{1}{\cos^2(x)}-1},
\end{equation}
although care must be taken as the square root can generate a different branch. This leads to the nested tangent representation 

\begin{equation}\label{eq:1}
\tan(x) = \sqrt{-1+{\lim_{n\to\infty}\Bigg\{-1\ldots+2\ldots\Big(-1+2\Big(-1+2\Big(1-\frac{x^2}{2^{2n+1}}\Big)^2\Big)^2\Big)^{2\ldots}\Bigg\}}^{-2}},
\end{equation}
and we can keep track of the branch due the sine function being positive or negative. Furthermore, the inverse is written as

\begin{equation}\label{eq:1}
\tan^{-1}(y) = \cos^{-1}\left(\frac{1}{\sqrt{1+y^2}}\right),
\end{equation}
which is represented by the inverse cosine function as

\begin{equation}\label{eq:1}
\tan^{-1}(y)= \lim_{n\to\infty} 2^{n}\sqrt{2\left(1-\sqrt{\frac{1}{2}\left(1+ \sqrt{\frac{1}{2}\left(1+\sqrt{\frac{1}{2}\left(1+\sqrt{\frac{1}{2}\left(1+\ldots \frac{1}{\sqrt{1+y^2}}\right)}\right)}\right)}\right)}\right)}.
\end{equation}

Finally, we remark on the exponential function. By inverting the nested logarithm function (70), we eventually obtain 

\begin{equation}\label{eq:1}
e^x=\cosh(x)+\sinh(x)=\cosh(x)+\sqrt{-1+\cosh^2(x)},
\end{equation}
hence basically, the nested cosine functions are fundamental in the sense of being able to invert them algebraically. It is not clear if the other trigonometric functions have an independent nested representation that can be inverted. We derived different variations which are coupled to the cosine function through different identities. Also, the exponential and logarithms function have different limit representations.  The most basic definition for the exponential is

\begin{equation}\label{eq:1}
e^x=\lim_{n \to \infty}\left(1+\frac{x}{n}\right)^n,
\end{equation}
which can be inverted to find the logarithm as
\begin{equation}\label{eq:1}
\log(y)=\lim_{n \to \infty}n\left(y^{\frac{1}{n}}-1\right)
\end{equation}
are the natural inverse limits.

\texttt{Email: art.kawalec@gmail.com}

\newpage
\section{Appendix A}

\textbf{Cosine function:}
\begin{equation}\label{eq:1}
\cos(x) = \lim_{n\to\infty}\Bigg\{-1\ldots+2\ldots\Big(-1+2\Big(-1+2\Big(1-\frac{x^2}{2^{2n+1}}\Big)^2\Big)^2\Big)^{2\ldots}\Bigg\}.
\nonumber
\end{equation}
\textbf{Inverse cosine function:}

\begin{equation}\label{eq:1}
\cos^{-1}(y)=\lim_{n\to\infty} 2^{n}\sqrt{2\left(1- \sqrt{\frac{1}{2}\left(1+\sqrt{\frac{1}{2}\left(1+\sqrt{\frac{1}{2}\left(1+\sqrt{\frac{1}{2}\left(1+\ldots y\right)}\right)}\right)}\right)}\right)}.
\nonumber
\end{equation}
\textbf{Sine function:}
\begin{equation}\label{eq:1}
\sin(x) = \sqrt{1-\lim_{n\to\infty}\Bigg\{-1\ldots+2\ldots\Big(-1+2\Big(-1+2\Big(1-\frac{x^2}{2^{2n+1}}\Big)^2\Big)^2\Big)^{2\ldots}\Bigg\}^2}.
\nonumber
\end{equation}
\textbf{Inverse sine function:}
\begin{equation}\label{eq:1}
\sin^{-1}(y)= \lim_{n\to\infty} 2^{n}\sqrt{2\left(1-\sqrt{\frac{1}{2}\left(1+ \sqrt{\frac{1}{2}\left(1+\sqrt{\frac{1}{2}\left(1+\sqrt{\frac{1}{2}\left(1+\ldots \sqrt{1-y^2}\right)}\right)}\right)}\right)}\right)}.
\nonumber
\end{equation}
\textbf{Tangent function:}
\begin{equation}\label{eq:1}
\tan(x) = \sqrt{-1+{\lim_{n\to\infty}\Bigg\{-1\ldots+2\ldots\Big(-1+2\Big(-1+2\Big(1-\frac{x^2}{2^{2n+1}}\Big)^2\Big)^2\Big)^{2\ldots}\Bigg\}}^{-2}}.
\nonumber
\end{equation}
\textbf{Inverse Tangent function:}
\begin{equation}\label{eq:1}
\tan^{-1}(y)= \lim_{n\to\infty} 2^{n}\sqrt{2\left(1-\sqrt{\frac{1}{2}\left(1+ \sqrt{\frac{1}{2}\left(1+\sqrt{\frac{1}{2}\left(1+\sqrt{\frac{1}{2}\left(1+\ldots \frac{1}{\sqrt{1+y^2}}\right)}\right)}\right)}\right)}\right)}.
\nonumber
\end{equation}
\textbf{Hyperbolic cosine function:}
\begin{equation}\label{eq:1}
\cosh(x) = \lim_{n\to\infty}\Bigg\{-1\ldots+2\ldots\Big(-1+2\Big(-1+2\Big(1+\frac{x^2}{2^{2n+1}}\Big)^2\Big)^2\Big)^{2\ldots}\Bigg\}.
\nonumber
\end{equation}
\textbf{Inverse hyperbolic cosine function:}
\begin{equation}\label{eq:1}
\cosh^{-1}(y)=\lim_{n\to\infty} 2^{n}\sqrt{2\left(-1+ \sqrt{\frac{1}{2}\left(1+\sqrt{\frac{1}{2}\left(1+\sqrt{\frac{1}{2}\left(1+\sqrt{\frac{1}{2}\left(1+\ldots y\right)}\right)}\right)}\right)}\right)}.
\nonumber
\end{equation}
\textbf{Hyperbolic sine function:}
\begin{equation}\label{eq:1}
\sinh(x) = \sqrt{-1+\lim_{n\to\infty}\Bigg\{-1\ldots+2\ldots\Big(-1+2\Big(-1+2\Big(1+\frac{x^2}{2^{2n+1}}\Big)^2\Big)^2\Big)^{2\ldots}\Bigg\}^2}.
\nonumber
\end{equation}
\textbf{Inverse hyperbolic sine function:}
\begin{equation}\label{eq:1}
\sinh^{-1}(y)= \lim_{n\to\infty} 2^{n}\sqrt{2\left(-1+\sqrt{\frac{1}{2}\left(1+ \sqrt{\frac{1}{2}\left(1+\sqrt{\frac{1}{2}\left(1+\sqrt{\frac{1}{2}\left(1+\ldots \sqrt{1+y^2}\right)}\right)}\right)}\right)}\right)}.
\nonumber
\end{equation}
\textbf{Hyperbolic tangent function:}
\begin{equation}\label{eq:1}
\tanh(x) = \sqrt{1-{\lim_{n\to\infty}\Bigg\{-1\ldots+2\ldots\Big(-1+2\Big(-1+2\Big(1+\frac{x^2}{2^{2n+1}}\Big)^2\Big)^2\Big)^{2\ldots}\Bigg\}}^{-2}}.
\nonumber
\end{equation}
\textbf{Inverse hyperbolic tangent function:}
\begin{equation}\label{eq:1}
\tanh^{-1}(y)= \lim_{n\to\infty} 2^{n}\sqrt{2\left(-1+\sqrt{\frac{1}{2}\left(1+ \sqrt{\frac{1}{2}\left(1+\sqrt{\frac{1}{2}\left(1+\sqrt{\frac{1}{2}\left(1+\ldots \frac{1}{\sqrt{1-y^2}}\right)}\right)}\right)}\right)}\right)}.
\nonumber
\end{equation}
\textbf{Exponential function:}
\begin{equation}\label{eq:1}
e^x=\cosh(x)+\sinh(x)=\cosh(x)+\sqrt{-1+\cosh^2(x)}
\nonumber
\end{equation}
\begin{equation}\label{eq:1}
e^x=\lim_{n \to \infty}\left(1+\frac{x}{n}\right)^n.
\nonumber
\end{equation}
\textbf{Logarithm function:}
\begin{equation}\label{eq:1}
\log(y) =\lim_{n\to\infty} 2^{n}\sqrt{2\left(-1+\sqrt{\frac{1}{2}\left(1+ \sqrt{\frac{1}{2}\left(1+\sqrt{\frac{1}{2}\left(1+\sqrt{\frac{1}{2}\left(1+\ldots \frac{y+y^{-1}}{2}\right)}\right)}\right)}\right)}\right)}
\nonumber
\end{equation}
\begin{equation}\label{eq:1}
\log(y)=\lim_{n \to \infty}n\left(y^{\frac{1}{n}}-1\right).
\nonumber
\end{equation}

\end{document}